\documentclass{ewic}

\usepackage{amssymb, amsmath}

\input turing.paper.defs

\begin{document}

\header{The Incomputable Alan Turing}

\footer{\em{Turing 2004: A celebration of his life and 
achievements}}

\title{The Incomputable Alan Turing}

\author{S.\ Barry Cooper\\
School of Mathematics, University of Leeds\\
Leeds, LS2 9JT, UK\\
www.maths.leeds.ac.uk/\,\~{}pmt6sbc\\
\email{s.b.cooper@leeds.ac.uk}}

\begin{abstract}
The last century saw dramatic challenges to the Laplacian
predictability which had underpinned scientific research for
around 300 years. Basic to this was Alan Turing's 1936
discovery (along with Alonzo Church) of the existence of
unsolvable problems. This paper focuses on incomputability as
a powerful theme in Turing's work and personal life, and
examines its role in his evolving concept of machine
intelligence. It also traces some of the ways in which
important new developments are anticipated by Turing's
ideas in logic.
\end{abstract}


\section{A beginning \dots and an end?}

Alan Turing's life and work had a peculiarly intimate 
relationship with each other. 
  Thanks to Andrew Hodges' insightful biography \cite{hodges}, 
it is possible to trace many of the details of this, and to 
place it within the context of  twentieth-century
 erosion of certainty.

{\rightskip=14pt \leftskip=14pt \hskip 0mm{\sl 
``I am making a collection of experiments in the order I mean to 
do them in. I always seem to want to make things from the thing that 
is commonest in nature and with the least waste in energy."}
\vskip 1mm}
 
So wrote (quoted by Hodges, p.19) the young Alan Turing in a  
 letter home from boarding school in March 1925, setting out 
the essentials of a lifelong agenda --- one whose 
shadowy outcome, with the very nature of predictability 
in question, could not have been foreseen.  
Here (wrenched out of context) 
is the closing paragraph (p.527) of  
Hodges' book: 

{\rightskip=14pt \leftskip=14pt \hskip 0mm {\sl 
``With so few messages from the unseen mind to work on, 
[Alan Turing's] inner code remains unbroken. According to 
his imitation principle, it is quite meaningless to speculate upon 
his unspoken thoughts. {\it Wovon man nicht 
sprechen kann, dar\"uber muss man schweigen.\/}\footnote{\,Usually 
translated as:  
``What we cannot speak about we must pass over in silence." 
--- being the final sentence of Ludwig Wittgenstein's
\textit{Tractatus Logico-Philosophicus\/}. }  
But Alan 
Turing could not possess the philosopher's detachment from life.
 It was, as the computer might put it, the unspeakable that left 
him speechless."}\vskip 1mm} 

What had happened by the end was that Turing had played a key role in 
 pushing forwards the boundaries of machine intelligence, 
had struggled with newly discovered theoretical and practical obstacles to  
 computability, and humanly 
confronted the parallel realities of a world out of control. And just as it 
was hard for Turing to make sense of what he had played a role in 
revealing, we have to be cautious, fifty years later, in interpreting and 
simplifying 
Turing's own views.

Some of the themes addressed in this article will be:   Turing's investment in 
computability,   incomputability,  
and in their practical dimensions;
his work in logic, and the ever present real world connection; and 
the relevance today of the Turing approach 
to logic   --- and   its anticipation of current research themes.

\section{Confronting the incomputability barrier --- the Alan Turing way}
Incomputability  has been with us for 
much longer than has the precise sense the term today derives from 
computability theory. Who has not made plans based on best calculations, 
only to be surprised and thwarted by the unpredictability of the real world?  
Even the word itself has been with us for around four hundred years. The 
{\sl New Oxford Dictionary of English\/} (1998) gives the meaning 
``unable to be calculated or 
estimated, {\sl origin\/} early 17th cent." 
In everyday usage, incomputability is a 
barrier to human activity --- something we all experience 
in everyday life. \smallskip

At this level, the commonest coping strategy is a cultural one --- 
imitation. One adopts 
strategies learned from  
copying what works for other people. A more individual, and 
necessarily rational, one is that 
of reductionism and thought experiments.  
 Here, complex phenomena are broken down into basic ones, 
the aim being solutions we recognise 
as being  scientific. We can describe this as `computing with awareness'.

 For Alan Turing, strategy I was not a favoured option\,! 
For him, the experience of incomputability   was 
an intimate one,  his scientific and  
personal life inextricably intertwined. 
To this extent his mind became a ready-made laboratory 
for exploring the computable\,/\,incomputable 
interface. Here is how his older brother John (Hodges, p.33) 
describes, around 1928, Alan's brave approach to making sense of the 
world:

{\rightskip=14pt \leftskip=14pt \hskip 0mm {\sl 
``You could take a safe bet that if you ventured on some 
self-evident proposition, as for example that the earth was 
round, Alan would produce a great deal of 
incontrovertible evidence to prove that it was almost 
certainly flat, ovular, or much the same shape as a Siamese 
cat which had been boiled for fifteen minutes 
at a temperature of one thousand degrees Centigrade."}\vskip 1mm}

\section{ \dots and finding this barrier all-too-real}
 On Thursday  13 February 1930, Turing's school friend 
Christopher Morcom died. This unexpected  blow to the young 
Alan Turing (see Hodges) was a seismic shattering of his 
personal Laplacian 
 universe. Alan tried to make sense of what had happened, 
 writing to his best friend's mother, Mrs.\ Morcom, around 1932   
(quoted in  Hodges, 
 pp.\ 63--64):

{\rightskip=14pt \leftskip=14pt \hskip 0mm {\sl
``It used to be supposed in Science that if 
everything was known about the Universe at any particular moment 
then we can predict what it will be through all the future. \dots 
 More modern science however has come to the 
conclusion that when we are dealing with atoms and electrons we 
are quite unable to know the exact state of them \dots The 
conception then of being able to know the exact state of the 
universe then really must break down on the small scale. This means 
then that the theory which held that as eclipses etc.\ are 
predestined so were all our actions breaks down too. We 
have a will which is able to determine the action of \dots  
atoms\dots in \dots the brain,\dots {\hskip 10mm}
\dots matter is meaningless in the absence of spirit \dots 
as regards the actual connection between spirit and body I consider 
that the body by reason of being a living body can 
`attract' and hold on to a `spirit' \dots The 
body provides something for the spirit to look 
after and use."}\vskip 1mm} 

It is worth looking more closely 
at this for what it tells us about Turing's 
thinking at this time. One is struck by a fluidity  of thought, anticipating 
current concerns (and confusions!) about computability. There 
are tentative steps towards a view of the Universe as {\sl information\/}. 
There is the
awareness of, and interest in, the then very recent  discovery of quantum uncertainty.
For the young Turing (familiar   with Eddington \cite{eddington}, for example)  there is
a very  contemporary resonance 
between the  personal and 
the scientific. {\sl Human spirit\/} becomes admissible, 
 embodying  
{\sl special\/} information content, from which it gets purpose and 
persistence.  And there are those last three wonderful clauses, 
full of a very contemporary sense of the mystery and significance 
of the relationship between different levels of complexity, 
and even a premonition of the underlying phenomenon of 
mathematical definability
detected  in recent explanations. 
 And of the {\sl necessity\/} of incomputability. 
Notice that at this point the mind is not seen as 
mechanical.

Is the above quotation little more than the unformed and naive
writings of a young man  getting through an emotionally stressful
part of his life?  One can so dismiss them. But only three or four
years later,  Turing was working on one of his  most admired and
influential  scientific papers. 

\section{ That remarkable 1936 paper} \vskip 1mm

Turing \cite{turing36} did not actually appear 
until January 1937. Typically, Turing starts with a 
more basic question than that  asked by other authors. Not ``What is a computable
function?". But (p.249): 

{\rightskip=14pt \leftskip=14pt \hskip 0mm {\sl
`The real question at issue is ``What are the possible processes 
which can be carried out in computing a [real] number?" '}\vskip 1mm} 

For a full commentary on this paper and its historical context, 
see Robin Gandy \cite{gandy} --- this is essential reading, of course, 
from Turing's ex-student and friend. Here is a familiar 
picture (taken from Cooper \cite{cooper04}) of Turing's hardware,  
all computational subtlety residing in the program.

\vskip 4mm 

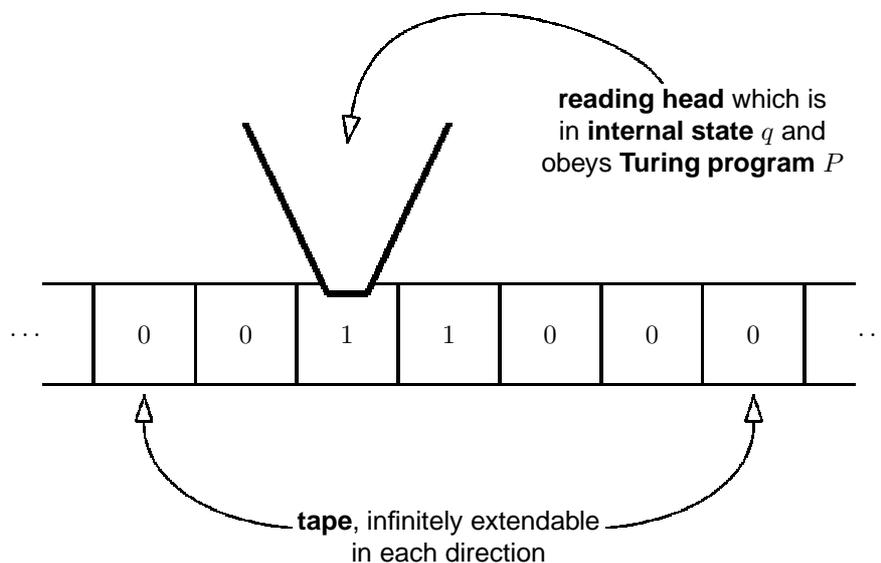
\Draw\Scale(1.9,1.9)

\LineAt(-10,0,45,0)
\LineAt(55,0,150,0)
\LineAt(-10,-20,150,-20)
\LineAt(0,0,0,-20)
\LineAt(20,0,20,-20)
\LineAt(40,0,40,-20)
\LineAt(60,0,60,-20)
\LineAt(80,0,80,-20)
\LineAt(100,0,100,-20)
\LineAt(120,0,120,-20)
\LineAt(140,0,140,-20)

\PenSize(.75mm)
\LineAt(30,32,46,-2)
\LineAt(46,-2,54,-2)
\LineAt(70,32,54,-2)

\MoveTo(10,-10)\Text(--$0$--)

\MoveTo(30,-10)\Text(--$0$--)

\MoveTo(50,-10)\Text(--$1$--)

\MoveTo(70,-10)\Text(--$1$--)

\MoveTo(90,-10)\Text(--$0$--)

\MoveTo(110,-10)\Text(--$0$--)

\MoveTo(130,-10)\Text(--$0$--)

\MoveTo(154,-10)\Text(--$\dots$--)
\MoveTo(-13,-10)\Text(--$\dots$--)

\MoveTo(118,30)\Node(1)(--{{\bf reading head} which is}
{in {\bf internal state} $q$ and}
{obeys {\bf Turing program} $P$}--)
\MoveTo(50,28)\FcNode(2)

\MoveTo(70,-50)
\Node(3)(--{{\bf tape}, infinitely extendable} 
{in each direction}--)
\MoveTo(10,-22)\FcNode(4)
\MoveTo(130,-22)\FcNode(5)

\ArrowHeads(1)
\PenSize(0.5pt)

\ArrowSpec(H)
\CurvedEdgeSpec(22,0.2,-66,0.2)
\CurvedEdge(3,4)

\CurvedEdgeSpec(-22,0.2,66,0.2)
\CurvedEdge(3,5)

\CurvedEdgeSpec(-60,0.2,80,0.5)

\CurvedEdge(1,2)

\EndDraw\vskip 4mm

The Turing machine may not have been the first  
general notion of computability --- Post even had a 
{\sl machine\/} notion of computability (not published) which was very 
similar to Turing's. 
 But it was the most convincingly formulated 
and presented. It won over 
G\"odel (see John Dawson \cite{dawson}, pp.\ 101--102). 
G\"odel's long-time friend Hao Wang provides (\cite{wang}, 
p.96), the final word 
on Kurt G\"odel and Church's thesis:

{\rightskip=14pt \leftskip=14pt \hskip 0mm {\sl 
``Over the years G habitually credited A.\ M.\ Turing's 
paper of 1936 as the definitive work in capturing the intuitive 
concept [of computability], and did not 
mention Church or E.\ Post in this connection. He must have felt 
that Turing was the only one who gave persuasive arguments to show the 
adequacy of the precise concept \hskip 3mm  \dots\hskip 4mm In particular, 
he had probably been aware of the arguments offered by 
Church for his `thesis' and decided that they were inadequate. 
It is clear that G and Turing (1912--1954) had great 
admiration for each other, \hskip 3mm \dots"}\vskip 1mm} 

Over the years, of course, the proved equivalences between   
the diverse models of computability gave plenty of evidence 
that the Church-Turing Thesis (as it is now called) is 
essentially correct. But it is the persuasiveness of 
Turing's original argument which still has more powerful  
intuitive content for us.

Turing's paper is characteristically basic in other ways.  
For instance, it also takes account of the scientist's need to describe 
the world in terms of {\it real numbers}. 
Turing's `computable real numbers' embodied concessions that the 
soon-to-become dominant recursive function theoretic 
framework clearly did not. 

\section{The Universal Machine and the Computer}

Not everything in the paper was so accessible. The student 
who imagines it was simply Turing's 
 invention of the Turing machine which anticipated the 
computer is missing the main point of his paper. 
Turing's idea was to G\"odel number (or code) programs, 
while building a {\sl universal\/} Turing machine $U$ capable of 
{\sl decoding\/} and {\sl implementing\/} 
any coded 
program.

As Martin Davis \cite{davis} points out --- 
 at a stroke, Turing anticipated:

\It {\sl Interpretive programs\/} --- in that $U$ unscrambles codes 
into implementable information (quintuples),

\It The {\sl Stored Program Computer\/} --- whereby the distinction 
between program and data evaporates,

\It The {\sl versatility\/} of today's computers, with {\sl
hardware\/}  solutions 
replaced by equivalent software.

It is clear from Turing's paper that  his universal 
machine is only universal in a very specific sense.   Turing's machines --- 
and the $U$ which simulates them --- 
were  framed in terms of the computing capabilities of a human clerk.  
Those of  machines in general  remained to be fully pinned 
down. Turing's later involvement with physical computability 
was to bring to the surface some contemporary themes 
still occupying  
those seeking to mathematically model how the real world 
computes.  

What is also 
 clear,  is that the universal Turing machine provides an 
accurate   
 basic model of how the computer on our 
desk computes.  From the perspective of the pure mathematician, Turing {\sl
did\/} invent the computer. You can write down the program for 
$U$ and implement it in a full range of contexts, admittedly 
in too inefficient a way for it to have any practical use\,! 
Turing's name and reputation still has power to enrage 
   the engineer 
concerned with implementing that  abstract idea. Turing's own 
involvement in the building of those first computers just 
makes his role even more irritating to those 
who neither recognise, nor want to recognise, 
the mathematical structures which constrain their daily lives.

Turing's best-known early influence on computer 
design was through John von Neumann, particularly via the latter's  
  EDVAC report (1945). Although 
von Neumann does not specifically 
acknowledge the significance of Turing there, he later 
gave proper credit in his Hixon Symposium lecture, in 
September 1948 (see \cite{neumann}). Understanding of Turing's key 
role in the design
of early  computing machinery, such as the `Colossi' (circa 1944), during the 
Second World War, is to a great extent lost through Churchill's 
deliberate obliteration of physical evidence of the Bletchley Park 
story. 

 But relevant for us today, as Robin Gandy \cite{gandy01} tells us (p.17) in his 
Preface to the 1936--37 papers in the {\sl Collected 
Works of A.M. Turing: Mathematical Logic\/}:

{\rightskip=14pt \leftskip=14pt \hskip 0mm {\sl
 ``Turing machines \dots still provide the standard setting for the 
definition of the complexity of computation in terms of bounds on time 
and space;\endgraf
$\bullet$ together with the neural nets of McCulloch and Pitts they provided 
the foundations of the theory of automata;\endgraf
$\bullet$ together with the generated sets of Post [1943] they provided 
the foundation for the theory of formal grammars."}\vskip 1mm}

\section{Hilbert's Programme, and Natural Examples of Incomputable Objects}

Turing's  paper was not immediately useful, of course. Its 
practical benefits lay still some years in the future. 
But it did relate to a grand scientific enterprise traceable 
back to Newton and beyond, that of  
capturing the algorithmic content of the natural world. As 
Albert Einstein
\cite{einstein} said (p.54): 

{\rightskip=14pt \leftskip=14pt \hskip 0mm {\sl
``When we say that we understand a group of 
natural phenomena, we mean that we have found a 
constructive theory which embraces them."}\vskip 1mm}

And for David Hilbert, this included    
mathematics and its epistemology.  
 {\sl Hilbert's Programme\/} (developed roughly over the 
period 1904--1928) was (to put it simply) to 
capture mathematics 
in complete, consistent theories. Hilbert's `Entscheidungsproblem' 
 asking (again in simple terms) for  an algorithm 
for deciding if a given sentence 
 is logically valid or not --- was an important 
aspect of the programme. We can see it as 
 growing out of Hilbert's view that:

{\rightskip=14pt \leftskip=14pt \hskip 0mm {\sl
 ``For the mathematician there is no
Ignorabimus, and,  in my opinion, 
not at all for natural science either. \dots The true reason 
why [no one] has succeeded in finding an unsolvable problem is, 
in my opinion, that there \underbar{is no} unsolvable problem. 
In contrast to the foolish Ignorabimus, our credo avers:\endgraf
We must know,\endgraf
We shall know."}\vskip 1mm}

Already, the day before this opening 
address to the Society of German Scientists 
and Physicians
in K\"onigsberg,  September 1930,  G\"odel had announced, at a different 
meeting in another part of the same city, the 
incompleteness result that was to have such 
negative impact on Hilbert's programme. And in 1936, 
it turned out that not only could one not fully capture 
mathematical truth in axiomatic theories but, worse,  
that which was captured was, in a sense, not captured 
at all\,!

 Church showed no such algorithm for the Entscheidungsproblem 
existed. And so did 
Turing --- essentially showing:

{
(1) The set of inputs $n\in\Bbb N$ on which $U$ halts is computably 
enumerable but not computable, and hence:\endgraf
(2) The logically valid sentences form an incomputable  c.e.\ set.
}

This was just the beginning. Soon a whole repertoire of incomputable c.e.\ sets
appeared,  and any reasonably rich theory turned out to be 
undecidable. Even questions of interest to `real' mathematicians 
were caught in the net. Turing himself played a part 
in the proving of the unsolvability of
the  word problem for groups (Post/ Markov 1947, Turing 1950, 
Novikov/ Boone 1955). 

There are two footnotes to the story of Turing's 1936 paper.

 We should remember there was a second, 1937 paper, 
of course. Turing was prone to technical 
mistakes, and these needed correcting. Gandy (\cite{gandy01}, Preface to the 1936--37
papers, p.9) 
nicely draws together the personal and the  scientific again, always
important with Turing (my  underlining):

{\rightskip=14pt \leftskip=14pt \hskip 0mm {\sl
``Whether calculating mentally or 
with pencil and paper, Turing was methodical only by fits and 
starts, and often made mistakes. [When I came to know him later 
the phrase `What's a factor of two between friends?' had 
become a catchword.]  But he understood 
very well what it meant to be totally methodical. Indeed 
\underbar{an acceptance} --- sometimes ready, sometimes reluctant --- 
\underbar{of the dichotom}y\underbar{ between the}\underbar{ clearl}y
p\underbar{erceived ideal and the 
confused actualit}y\underbar{ was fundamental in Turin}g\underbar{'s thought}."}\vskip
1mm}

And it was Alonzo Church   
who first  proved Church's Theorem (as the negative solution to the 
Entscheidungsproblem came to be called). This time it is Hodges 
(\cite{hodges}, p.114) 
who finds the scientific pearl in the Turing 
personal oyster:

{\rightskip=14pt \leftskip=14pt \hskip 0mm {\sl
``If he had been a more conventional 
worker, he would not have attacked the Hilbert problem without 
having read up all of the available literature, including 
Church's work. He then might not have been pre-empted --- 
but then, he might never have created the new idea of the logical 
machine, with its simulation of `states of mind', which not 
only closed the Hilbert problem but opened up quite 
new questions."}\vskip
1mm}

\section{But incomputability into the long grass \dots}

Now we come to one of the oddest and most 
interesting things about the 
 the 1936 discoveries and their repercussions. The 
emergence of incomputability soon after that 
of quantum unpredictability, with its iconoclastic 
import for old certainties, was hardly acknowledged outside 
the hermetic confines of pure mathematics. At the same time, the 
universal Turing machine, a mere by-product of that 
unexpected glimpse beyond Laplacian predictability, 
became the basis of an informational  
revolution built on filling out the computable (see, 
for example, Chapter 7 of \cite{davis}). 

What happened, and why? Was what happened the natural order of things, 
or can we see some shape and development underlying these events?

A first observation is that 1936 not only saw a new clarity about 
what `incomputability' really
is, but was accompanied by the emergence of a conceptual framework 
which actually took  us {\sl away\/} from 
the real world and the uncertainties 
facing working scientists. Sometimes, being precise about 
notions has negative as well as positive consequences. It may make  
them technically more useful, while disrupting their 
more general usefulness. 
In the years following 1936 we saw  
the birth of {\sl Recursion Theory\/}, and
new notions inimical  to vague intuitions.  

This involved a dual defeat. Not just that of 
the outsider trying --- not very hard it often seemed\,! --- 
to understand what recursion theorists were doing, and why, and failing. But 
of the specialist looking out, usually when seeking funding.   
One saw a growing emphasis on purely
{\sl mathematical issues\/},  extending to logic in general. 
The technical work was often accompanied by an overpowering 
feeling of reality and subjective significance --- but 
guided by  an 
introjected real world rather than an actual one. A number of those with 
a professional interest in incomputability (and in 
large cardinals, for that matter) were even convinced   
finitists in relation to the material universe.  Not for them  
far-reaching speculations about how the universe, or the human mind, 
transcends the Turing-barrier, at least not in public. 
On the other hand, some wonderful mathematics resulted (see 
Odifreddi's monumental work \cite{od89}, \cite{od99}), and
this  isolated development, one suspects,  
was entirely necessary to the future relevance of the subject. 
Only now are new connections being made (see, for example,  
\cite{CoOd03}), based on the 
very hardest technical work.

At the same time, there seemed to be a 
 growing belief that mathematics --- and science 
in general ---  could carry on much as before  without 
ever bumping into incomputable objects.  
Laplacian  predictability was not just a useful model, 
it was a lasting mind-set amongst working scientists.  

There were also mathematical developments which 
suggested incomputability might be just as containable 
as sub-atomic unpredictability had become within quantum theory.  
John Myhill \cite{myhill} showed that all the
unsolvable problems  discovered in the 1930s were 
essentially {\sl the same\/}. They were just notational 
variants of what he called a {\sl creative\/} set. 
One could artificially manufacture other incomputable 
objects, but there was no evidence of their existence outside of 
academic papers. 
Perhaps 
incomputability was not as pervasive and naturally occuring as it had at 
first seemed \dots

Anyway, there seemed little point in worrying about the recursion theorists' negative
news, when  there was so much computability to explore.   
The  richness of the {\sl
computable\/} universe was increasingly revealed, and provided 
more than enough work for a whole army of researchers. 
And there was even a mathematical dimension to 
this.  {\sl Reverse
mathematics\/} emerged as an attempt to rescue  Hilbert's programme 
--- with some success, it seemed to many. When leading 
proponent Harvey Friedman challenged the recursion theorists to 
produce {\sl just one\/}  natural example of an 
incomputable 
object (see, for example \cite{friedman}), he chose his own 
rules --- a
peculiarly  mathematical definition of ``natural" --- which left 
them with nowhere to turn.

\section{The Search for \underbar{Reall}y Natural Examples}

Back in the real world, there is a definition of ``natural" 
understood by everyone. Not a metaphorical one, but one 
of quotidian usefulness. The {\sl 
New Oxford Dictionary of English\/}\, 
(1998 edition), gives it to us as:

{\rightskip=14pt \leftskip=14pt \hskip 0mm
  {``{\bf natural} --- {\sl existing in or caused by 
 nature\/}"}\vskip
1mm}

Numbers --- meaning the {\sl natural\/} numbers --- 
are generally thought to qualify here. In fact, the whole of 
school arithmetic does. In 1972, Yuri Matiyasevich put the 
finishing touches to a remarkable proof, many years 
in the making, that the sort of arithmetical questions 
a school student might ask lead directly to a rich 
diversity of incomputable sets --- a diversity which 
went far beyond those original examples from the 1930s.
Once again, it was David Hilbert being audaciously wrong which 
provided the impetus, and drew public attention to the 
result. The negative solution to Hilbert's Tenth Problem,  
collectively due to Martin Davis, Matiyasevich himself, 
the philosopher Hilary Putnam, and, crucially, 
Julia Robinson, demonstrated that  all the 
incomputable sets which could be artificially enumerated 
by clever logical techniques were not artificial at all --- 
they arose {\sl naturally\/} as solution sets of 
the familiar diophantine equations. 
  It is {\sl everyday\/} mathematics leads us unavoidably to 
incomputable objects. 

For the working scientist, naturalness is captured within a rather 
different kind of mathematics, full of real numbers and differential 
equations. Marian Pour-El was one of the first to look 
for appropriate notions and evidence of incomputability 
arising from this context.  One of the most widely-known 
and enduring discoveries, due to 
 Pour-El and Richards, was a differential equation with computable 
boundary conditions leading to incomputable solutions. 
More recently, high-profile names (such as Roger Penrose, Steve Smale) 
have been associated with investigations of the computability of 
well-known mathematical objects with known 
connections with complexity in nature, such as the Mandelbrot 
and Julia sets. 

The closer one gets to physical reality, the more potentially 
persuasive, while at the same time more speculative and 
ultimately elusive, the examples become. Leaving aside for the moment the workings 
of the human brain, the most scrutinised and puzzling part 
of the physical world is the sub-atomic level. 
 The predictive incompleteness of quantum theory refuses 
to go away, giving rise to different `interpretations'  which 
leave us a long way from characterising the algorithmic 
content of the events it seeks to describe. The 
quantum process which seems to escape the predictive net 
most radically, the most promising avatar of incomputability, 
is by-passed by current quantum computational models --- 
despite recent claims (see 
Tien Kieu \cite{kieu}). This is Andrew Hodges' comment on the situation 
(taken from his article {\sl What would Alan Turing have done 
after 1954?\/}, in Teuscher \cite{teuscher}):

{\rightskip=14pt \leftskip=14pt \hskip 0mm {\sl
``Von Neumann's axioms distinguished the {\bf U} 
(unitary evolution) and {\bf R} (reduction) rules of quantum mechanics. 
Now, quantum computing so far (in the work of Feynman, Deutsch, 
Shor, etc\.) is based on the 
{\bf U} process and so computable. It has not made serious use of 
the {\bf R} process: the unpredictable element 
that comes in with reduction, measurement, 
or collapse of the wave function."}\vskip
1mm}

Observable signs of incomputability in nature are not so obvious at the 
classical level - we are more embroiled and cannot so easily 
objectify what we are part of. As early as 1970, Georg Kreisel was not deterred, 
and in a
footnote to \cite{kreisel} (p.143)  went so far as to propose the possibility of a   
 collision problem related to 
the 3-body problem  which might  
give ``an analog computation of a 
non-recursive function (by repeating collision 
experiments sufficiently often)". 

This conjecture has come to seem less outrageous as 
people from various backgrounds have come up with even more 
basic proposals. One of these comes out of recent progress 
on the  Painlev\'e Ê Problem, from 1897, asking whether  
 noncollision singularities exist for the $N$-body problem 
for $N\geq 4$.  For $N\geq 5$, Jeff Xia in 1988 showed the answer is ``Yes" 
  (see  Saari and Xia \cite{saari}). And actual, or even uncompleted, 
infinities in nature are what open the door to incomputability. 
There are a number of other recent examples.

\section{The Incomputable -- so near but so far\dots}

Turing himself returned to the topic of incomputability in  
  what is usually thought of as a rather opaque 1939 paper. 

The background situation established in the years 1931--36 puzzled many people, 
including Turing himself. They knew that no Turing machine could 
prove all the true sentences  of arithmetic.  
(Turing's work had enabled a precise notion
of  {\it formal system\/} to apply  G\"odel's incompleteness 
procedure to.) 
 But, it seemed, a human observer could transcend what any 
{\sl given\/} such machine 
could prove. 

An immediate question was: Is there a mathematical 
analysis throwing light on the {\sl apparent\/} ability of the human mind 
to transcend the mechanical --- hopefully based on   
an extended constructivism? G\"odel, in setting out to 
verify elements of Hilbert's programme, had found his 
incompleteness theorem. Now Turing, in trying to 
place the mathematician's thought processes within 
a constructive framework, was to reveal  
a similarly surprising scenario. 
Andrew Hodges (\cite{hodges}, 
p.137) refers to Turing's conversations with on his way back to America 
in 1937:

{\rightskip=14pt \leftskip=14pt \hskip 0mm {\sl
``\dots now [Turing] gave the impression that he had long been 
happy with the Russellian view, that at some level the world 
must evolve in a mechanistic way. \dots Symbolically, 
the Research fountain pen that Mrs Morcom had given him 
in 1932 was lost on the voyage."}\vskip
1mm}

Turing was visiting the Institute of Advanced Study in Princeton, 
and writing a PhD thesis under Alonzo Church --- one of the few 
 celebrated mathematicians he had hoped to meet there (such as 
G\"odel and von Neumann) who were actually in residence during his 
time there. It was this thesis which would provide the material 
for the 1939 paper \cite{turing39}. 

Roughly speaking, the 
key idea --- perhaps originating with 
Church --- was to use the  constructive ordinals $\mathcal O$
of Church and  Kleene (1937) to inductively extend theories
via       
  G\"odel-like unprovable sentences.   Turing was specially interested in 
bringing the true $\Pi^0_2$
sentences of arithmetic within this 
constructive framework. He observed that most mathematically interesting
problems,  such as  the Riemann hypothesis, are met at that level. 

His partial success was to get a hierarchy containing proofs 
for  all true $\Pi^0_1$ sentences of arithmetic. Even this 
outcome  
was not entirely satisfactory ---  
it turned out that  {\sl different\/} hierarchies of this kind can 
be complete (i.e., exhaustive), {\sl or\/} invariant (not vary  
with different notations for the same ordinal),  but not {\sl both}. 
Turing was interested in invariance, since he saw that as allowing one to 
unambiguously 
{\sl classify\/} problems according to their `depth' (i.e., level of 
ordinal notation). 
 Anyway,  the main disappointment was that the techniques did not seem to work for 
$\Pi^0_2$ sentences (this result awaited  Feferman and the use of 
 stronger reflection 
principles).

For many people, Turing's 1939 paper is best known for 
the first appearance there of {\sl oracle Turing machines\/}. 
 In investigating the 2-quantifier sentences, Turing sought a 
constructive derivation of a non-$\Pi^0_2$ problem --- 
and in so doing invented {\sl relativisation\/} (using 
 oracle machines) and, in essence, the {\sl Turing jump\/}\,!

This led to Post's wonderful 
1944 paper \cite{post44}, and 1948 short abstract \cite{post48}, 
which described some
first  far-reaching consequences (such as Post's theorem, the notion 
of Turing reducibility, the 
degrees of 
unsolvability), and  clarified  what is 
happening in the 1939 paper. 
 These three publications of Turing and Post  established the 
 still crucial theme of the 
{\sl interrelationship between computability and 
information content\/},  
 and established the  Turing universe  of 
algorithmically related reals  as the standard model for 
computationally complex environments.

\section{Computing the incomputable?}

Despite its technical complexities, made even less 
approachable by being based (presumably at Church's 
suggestion) on the lambda calculus model of 
computability, the 1939 paper 
has a characteristically mundane motivation. Turing claims to clarify 
here the 
relationship between `ingenuity' (subsumed 
within the ordinal logics) 
 and `intuition' (needed to identify good ordinal 
notations ---  an $\emptyset^{(\omega)}$ level 
intuition!) Turing clearly regards ingenuity as being 
what a clever Turing program is capable of, and 
intuition as something else. There is a clear implication 
that intuition is a feature of human mental processes, 
and to that extent Turing is certainly saying that 
his hierarchies have something to say about how the 
mathematician's mind transcends his own model of machine 
computability --- even if the results can be subsequently 
translated into proofs implementable by a Turing machine\,! 

But what about  `hypercomputation'? --- often associated  
nowadays with: {\sl Computing the incomputable using oracle Turing 
machines} (where this use may be implicit). The claim is (see, for example, Copeland 
\cite{copeland98}, \cite{copeland00}, 
Copeland and Proudfoot \cite{CopeProud99}), that Turing anticipated such a
possibility,  and that oracle Turing machines have been neglected 
in this respect until rediscovered in the late nineteen-nineties. 
Well, it takes a truly creative reading of the paper to 
give oracle machines such an explicit role in Turing's 
thinking about incomputability, at least at this point. That 
is not said with the dismissiveness of Martin Davis' 
recent articles (see, e.g., \cite{davis04}). It is certainly 
true --- see below --- that computers which interact 
 played an important part in Turing's 
later thinking, and that there was never any 
 evidence that Turing thought the human mind 
was limited to computing recursive functions.  
But no --- oracle machines, despite their 
undoubted relevance, are   
  tangential to Turing's thoughts here. 

Let us get back to what Turing (\cite{turing39}, 
  pp.134--5), actually says about 
the underlying meaning of his paper (my underlining): 

{\rightskip=14pt \leftskip=14pt \hskip 0mm {\sl
``Mathematical  reasoning may be regarded \dots as the exercise 
of a combination of \dots {\sl intuition\/} and 
{\sl ingenuity\/}. \dots In pre-G\"odel times it was thought by some that 
all the intuitive judgements of mathematics could be replaced by 
a finite number of \dots rules. The necessity for intuition would then be 
entirely eliminated. 
In our discussions, however, we have gone to the opposite 
extreme and eliminated not intuition but ingenuity, and this 
\underbar{in s}p\underbar{ite of the fact that} 
\underbar{our aim has been in much the same direction}."}\vskip
1mm}

So he is addressing the familiar mystery of how we 
often arrive at a mathematical result via what seems 
like a very unmechanical process, but then 
promptly retrieve from this a proof which is quite 
standard and communicable to other mathematicians. 
Another celebrated mathematician, well-known for his interest in the role of 
intuition in the mathematician's thinking, was Poincar\'e. 
A few years after Turing wrote the 
above passage, Jacques Hadamard \cite{hadamard} 
recounts how Poincar\'e got stuck 
on a problem (the content of which is not important):

{\rightskip=14pt \leftskip=14pt \hskip 0mm {\sl 
``At 
first Poincar\'e attacked [a problem] vainly 
for a fortnight, attempting to prove there could 
not be any such function \dots [quoting Poincar\'e:]}\vskip -2mm
{\sl Having reached Coutances, we entered an omnibus to go some 
place or other. At the moment when I put my foot on the step, 
the idea came to me, without anything in my former thoughts 
seeming to have paved the way for it \dots I did not verify the idea \dots 
I went on with a conversation already commenced, but I felt 
a perfect certainty. On my return to Caen, for conscience sake, 
I verified the result at my leisure."}\vskip
1mm}

 Who else but Turing would have attempted 
a mathematical explanation at that time? His 
argument is still not widely known, and its signficance 
certainly not understood, except by those at ease with both 
the mathematics and with thinking about the world 
in the sort of basic terms which came naturally 
to Turing. 

Emil Post (see \cite{post41}, p.55) was less subtle 
 in 1941, anticipating others 
such as Penrose:

{\rightskip=14pt \leftskip=14pt \hskip 0mm {\sl 
``\dots we may write\vskip 0mm
\centerline{\sl \underline{The 
Lo}g\underline{ical Process is Essentiall}y\underline{ Creative}}\vskip 0mm  
This conclusion, \dots makes of the mathematician much more than 
a kind of clever being who can do quickly what a \underline{machine} 
could do ultimately. We see that a \underline{machine} would 
never give a complete logic; for once the machine is made 
\underline{we} could prove a theorem it does not prove."}\vskip
1mm}

By the 1990s this sense of context is lost. By now 
it is no surprise to find this  
common amongst those  
of a recursion theoretic background.  Here is 
a typical quote from Stephen G. Simpson, taken from a 
widely circulated communication to the 
{\sl Foundations of Mathematics\/} (FOM) 
e-mail list, August 18 1998:

{\rightskip=14pt \leftskip=14pt \hskip 0mm {\sl 
``\dots Soare insists \dots that computability includes
relative computability, i.e. relative recursiveness as a means of
classifying non-recursive sets.  This terminology strikes me as
wrong-headed, as if one were to insist that biology includes the study
and classification of inanimate objects."}\vskip 1mm} 

Here is no longer any sense of engagement, beyond the purely 
mathematical, with  the  
complexities encountered at the interface between the computable and the 
incomputable. Robert Soare's call (e.g., \cite{soare}) for `recursion 
theory' to be returned to the computability-theoretic 
terminology of Turing is angrily rejected. 
The mathematics is deeply real to  
the writer quoted, but this reality is a hermetic one, the science 
unrecognisable to a latter-day Turing.  

This is not to say Turing's 1939 paper did not have 
a very important {\sl technical\/} influence within logic. This was 
very much the case in regard to proof theory.
For instance, Feferman et al were subsequently 
to replace Turing's consistency principles 
with much more powerful reflection principles, and  
(see \cite{feferman91}) 
much improve on Turing's main results. 

In fact, Feferman nicely clarified Turing's paper, 
including translating proofs out of the lambda-theoretic 
framework, in his 
1988 article \cite{feferman88} {\sl Turing in the Land of O(z)\/}. 
And on p.127 of that paper Feferman comments:

{\rightskip=14pt \leftskip=14pt \hskip 0mm {\sl 
`` \dots Turing anticipated \dots the classification by ordinals 
of the provably (total) recursive functions of various formal systems, 
obtained later by proof-theoretical work."}\vskip 1mm}

One can even find an anticipation of the 
 Paris-Harrington \cite{paris77} bringing of 
G\"odel's independence result for Peano arithmetic closer to home. 

This paper of Turing is a paper packed full of ideas, of course. In
section 10  (`The continuum hypothesis. A
digression')  Turing pursues a set-theoretic 
analogy, indicating how to replace $\omega_1$ by the
constructive  ordinals, and the subsets of $\omega$ by the computable reals --- 
so anticipating subsequent hierarchies of computable functions.

Later ---  Feferman's work \cite{feferman62} on `autonomous' ordinal logic  
(see Franz\'en \cite{franzen} for a readable account) puts 
Penrose's  brave but flawed speculations \cite{penrose} on computability 
and the mind in context.

\section{The Turing universe}

Turing's oracle machines give the real numbers an algorithmic 
infrastructure, which comprises the {\sl Turing universe\/}.  
Emil Post \cite{post48} gathered together the computably equivalent reals 
of this structure, and called the resulting ordering the 
{\sl degrees of unsolvability\/} --- later 
called the {\sl Turing degrees\/} (poor Post 
rarely got the credit) ---  and this has become the 
 mathematical context for the study of the Turing universe. 
And so arose one of the most 
forbidding and esoteric research topics in the whole of 
logic, or mathematics, for that matter.

It is only recently that attention has turned to Turing's universe of computably 
related  reals  as providing a model for  
scientific descriptions of a computationally complex 
real universe (see \cite{cooper99},  
\cite{CoOd03}, \cite{copeland98}, etc.) This is based on a growing appreciation 
of how algorithmic content brings with it  
 implicit infinities, and a science --- 
increasingly coming to terms with chaotic 
and non-local phenomena --- {\sl necessarily\/} 
framed in terms of reals rather than within some 
discrete or even finite mathematical model. 
There is a huge amount of research activity, much of it ad hoc in nature, 
concerned with the  
 computational significance of evolutionary and emergent 
form, and emergence in more specific contexts. And while 
 purists may claim this is `nothing to do with logic', there is 
undoubtedly 
 an urgent need here for the unifying and clarifying 
role of basic mathematical structures. (As is well-known, Turing himself 
was very interested in the emergence of form in nature, and wrote 
seminal papers on the topic --- for example \cite{turing52}.)

Of course, the persuasiveness of the full Turing model 
depends on the extent to which one 
sees around one  avatars of incomputability, and there 
is a strand of thought --- the hypercomputational,  
as Jack Copeland and others term it 
 --- concerned with contriving incomputability via 
explicitly physical versions of the Turing universe. 
What is common to both 
the computability-theoretic and 
hypercomputational strands is that both the emergence of incomputability, 
and the emergence of new relations in a universe 
which admits incomputability, are based on 
a better understanding of how the local and the global 
interact. Whatever the context, the key mathematical parallel here is that 
of {\sl definability\/} or {\sl invariance\/}, even if  
within rather different corresponding frameworks. 
This is not very explicit in building hypercomputational 
models, which enables Martin Davis and others to  
trivialise what is happening as being the use of oracles to  
shuffle around existing incomputability. 

Our renewed awareness of  
 real-world relevance puts our view of 
the pure mathematical theory of the Turing 
universe in a new light. 
An area of research which became known for 
its mathematical unlovelyness and 
forbidding pathology turns out to be the 
ideal counterpart to  real-world complexity. 

This was 
not at all clear during the recursion theoretic years. 
The difficulty of the area may have surrounded 
researchers of the 1960s --- pre-eminently Gerald Sacks --- 
with a vaguely heroic, even machismatic, aura. But as time 
went on this had become a double-edged weapon, and by the 
1990s almost no one was impressed by the length and incomprehensibility of 
groundbreaking new proofs. `Touching the Void' --- and having 
accidents --- was all very well for mountaineers  but,  
as the new century approached, 
mathematics  was very much about 
deliverables. 
At times the very value of research into 
relative computability was questioned. 
Images of such dissent stick in the memory:  
 Sacks himself, lecturing at Odifreddi's 
CIME summer school in Bressanone, Italy in 1979,   
illustrating his view of 
`Ordinary Recursion Theory' with a slide of  
 the Chinese masses in cultural revolution turmoil ---  
his metaphor for an activity obsessive, formless, pointless;  
or, ten years later, Robin Gandy's contribution to a discussion 
on the future of logic, at a conference in Varna, Bulgaria  --- communicating 
an impression of the  
structure
of the  Turing degrees via exaggeratedly 
 desperate  scribbles on a blackboard. 

By the end of the last century the subject would have 
been unrecognisable 
to a returning Alan Turing. 
We had a  Turing universe framed by failed mathematical 
ambitions, and isolated from its natural home, 
the complexity of the material world. 
And a reluctance by researchers to stray beyond purely technical 
questions, unlike Turing himself. 
The purely mathematical horizons 
of its brightest minds were 
focused on an esoteric enterprise 
centred around something called the 
 {\sl bi-interpretability conjecture\/}. The idea was to 
show that the Turing universe was really just   
  second-order arithmetic in disguise, a messy mathematical 
structure, but messy in a mathematically familiar, even respectable, way. A 
{\sl rigid\/} structure, its definable relations 
easily read off via their counterparts in true second-order 
arithmetic. Luckily for computability theorists and everyone else, 
this reductionism ran  
its course, and we still have  
a complex world  in need of 
the theoretical insights of the mathematical legatees of 
Alan Turing. 

Finally, I should 
mention that the conceptual leap which gave us the oracle computing machine 
had wider mathematical consequences. As Sol Feferman 
says in  \cite{feferman88} ({\sl Turing in the Land of O(z)\/}, p.127):

{\rightskip=14pt \leftskip=14pt \hskip 0mm {\sl 
``\dots eventually, the idea of 
transforming computability from an {\sl absolute\/} 
notion into a {\sl relative\/} notion would serve to open up the entire 
subject of generalized recursion theory."}\vskip 1mm}

And that of course is another story, which will have to be 
told elsewhere.

\section{Mind and the Church-Turing Thesis}

Turing himself was interested in the nature of physical 
computability right to the end. He may not have 
had in mind specific generalisations 
of the Church-Turing thesis --- see Odifreddi \cite{od89} 
for the best introduction to this topic --- but 
his preoccupation with the links between human and 
mechanical computability forced him to confront 
the underlying realities of physical computation. 
The following quote from Feferman's 
{\sl Turing in the Land of O(z)\/} (see \cite{feferman88}, 
pp.131--2), gives the usual impression people have 
of Turing's outlook:

{\rightskip=14pt \leftskip=14pt \hskip 0mm {\sl 
``Turing, as is well known, had a mechanistic conception 
of mind, and that conviction led him to have faith in the 
possibility of machines exhibiting intelligent 
behavior."}\vskip 1mm}

But these words have been carefully chosen, and if read as carefully 
do not prevent us following the often shifting balance in Turing's real world 
between logic and  
science. At one extreme we have Turing, founder of 
Artificial Intelligence and seminal influence on its 
methodology (e.g., via the Turing Test).  At the 
other we have his interest in
quantum theory, running through his writings for  Mrs Morcom, 
right up to his late postcards to Robin Gandy (see 
\cite{hodges}, p.\ 512).  And in between he considered  possibilities,  
coming out of 
his 1944-48 experiences  of the ACE (`Automatic 
Computing Engine') project,  such as 
machines which make mistakes. Here is Turing in his talk 
 to the London Mathematical Society, February 20, 1947  
 (quoted in Hodges, p.361):

{\rightskip=14pt \leftskip=14pt \hskip 0mm {\sl 
``\dots if a machine is expected to be infallible, it 
cannot also be intelligent. There are several theorems which say 
almost exactly that."}\vskip 1mm}

Turing also anticipated the importance now given to connectionist 
models of computation. See his discussion of `unorganised 
machines' in \cite{turing48}, and   
 Jack Copeland and Diane Proudfoot's article \cite{CopeProud98}  
``On Alan Turing's Anticipation of Connectionism".

A `mechanistic conception of mind' maybe, but no crude 
extension of the Church-Turing thesis in sight, even at a time 
when Turing had a huge personal 
investment in the development of 
computing machinery. What one sees in Turing's thinking 
is an attempt to relate computability-theoretic 
structure to physical reality in a very bold and 
basic way. Maybe he did not mention oracles after 1939, 
but he never ceased to emphasise the importance of 
{\sl context\/} for real-world computation. The 
following quote is from the same 
talk as the previous one (Hodges, p.361 again):

{\rightskip=14pt \leftskip=14pt \hskip 0mm {\sl 
``No man adds very much to the body 
of knowledge. Why should we expect more of a machine? 
Putting the same point differently, the machine must be allowed to 
have contact with human beings in order that it may 
adapt itself to their standards."}\vskip 1mm}

The mysteries Turing grappled with remain. 
To what extent is the logic of a Turing machine sufficient to capture 
the workings of a human brain?  
What is the nature of the mechanical in the physical world? And what 
relationship does this have to the  mind? 

He was to make one more major contribution to how we 
approach, and give meaning to, such difficult questions.

\section{Turing's seminal 1950 AI paper} 

Nowadays one's heart might sink on the reading of an introduction
to a   paper which proclaimed: 

{\rightskip=14pt \leftskip=14pt \hskip 0mm {\sl 
``I propose to consider the question, 
`\underbar{Can \, machines \, think} ?' "}\vskip 1mm}

We have become too used to inflated claims in this area, 
ones which turn out   
to repeat old arguments and, more often than not, old confusions. 
But Turing's 1950 paper \cite{turing50} on {\sl Computing machinery 
and intelligence\/}, in \textit{Mind\/}, both managed 
to say something new, and, as one expects from Turing, 
is as careful in formulating the question as it is 
in giving answers.     
 Despite what some people think, he did 
not ask ``Is the human mind a Turing machine", and still less 
did he answer, or even imply, ``yes" to such a question. 

What he did ask (\cite{turing50}, p.442) 
became the  basis for the standard approach to gauging the 
subsequent success of the AI enterprise:

{{\bf The Turing Test:}  {\sl `Are there 
imaginable digital 
computers which would do well in the imitation game ?'} 

Roughly speaking, the imitation game was described in terms of 
asking a human interrogator to distinguish between 
a human and a mechanical respondent, within a fixed span of time 
and under conditions which did not unfairly favour 
either of them. And in the paper, Turing does confess to a personal 
opinion. He makes (\cite{turing50}, p.442) a surprisingly  
limited claim for potential machine intelligence:

{\rightskip=14pt \leftskip=14pt \hskip 0mm {\sl 
``I believe that in about 
fifty years' time it will be possible to programme computers \dots 
to make them play the imitation game so well that an average 
interrogator will not have more than 70 per cent.\ chance of 
making the right identification after five minutes of questioning. The 
original question \dots I believe to be too meaningless to deserve 
discussion."}\vskip 1mm}

Whether or not Turing had a ``mechanistic conception 
of mind",  
his background 1938--39 experiences were not forgotten.
On the one hand he knew how G\"odel's theorem   
could create an illusion that people transcend computers, 
and this feeds into Turing's consideration of the 
  `The Mathematical Objection' to computers doing well 
at the imitation game. On the other hand, 
he would 
remember how his ordinal logics falsely suggest machines transcend 
the computable --- and their ultimate failure to `eliminate' intuition. 
The discussion is characteristically full of interest. 
For instance he makes a brave attempt 
at dealing with an objection based on `Continuity of the
Nervous  System'. And, yet again,  
he mentions  `Learning Machines'. Over the years,  
 oracles continued to feature implicitly  in Turing's thinking about 
machine intelligence. And at no time did he imply that the human 
mind could be modelled by a standard unaided Turing machine. 

While preparing the talk on which this article is based, I 
wrote to my old friend George\footnote{Piergiorgio} Odifreddi,
whose opinion on  this point I thought might be useful. This is an
extract  from what he wrote back 
to me on May 1st, 2004: 

{\rightskip=14pt \leftskip=14pt \hskip 0mm {\sl 
``\dots when reading turing's 1939 paper i DID have the
impression that he thought that by an oracle he meant a human
being, and thus that non-computable functions could be humanly
computable. the oracle device could be thought of as a
formalization of a human-machine interaction, in which the calls
to the oracle(s) would be a kind of human help received by the
machine. if this interpretation were correct, then it would mean
that turing did not accept the church-turing thesis that recursive
= human computability. quite the contrary, actually."}\vskip 1mm}

I am grateful to George for allowing me to
use what is, of course,  a hurriedly composed e-mail, although one
which I believe to be broadly  correct. The extensive literature
discussing mechanism, the mind, and 
 the work and opinions of G\"odel and Turing, 
does not get us much closer to a definitive view. Two recommendable recent  
contributions treading similar ground with differing 
outcomes are Shapiro \cite{shapiro98} and Piccinini \cite{piccinini03}. 

\section{The Turing Renaissance}

What is currently so exciting is that the sorts of questions 
which preoccupied Turing, and the very basic 
extra-disciplinary thinking which he 
brought to the area, are being revisited and renewed by researchers 
from quite diverse backgrounds. What we are 
seeing is an emergent coming together of logicians, 
computer scientists, theoretical physicists, people 
from the life sciences, and the humanities and beyond, 
around an  intellectually 
coherent set of computability-related problems.  The recurring and 
closely linked themes here are the relationship between the local and the global, 
the nature of the physical world, and within that the human mind, as a computing
instrument, and  our expanding concept of what may be practically computable. 

The specific form in which these 
themes become manifest are quite varied. 
For some there is a direct interest in incomputability in Nature, 
such as that coming out of the n-body problem or  
quantum phenomena. For others it is through addressing  
 problems computing with reals and with scientific computing. 
 The possibility of computations `beyond the Turing barrier' 
leads to the study of    
analog computers, while theoretical 
models of hypercomputation figure in heated 
cross-disciplinary controversies. 
There is also intensive research going on into 
a number of practical models of natural computing, 
which present new paradigms of computing whose 
exact content is as yet not fully understood. 
In many scientific areas the emergence of form 
is deeply puzzling, and there is a need for mathematical models.  
Turing played an anticipatory role 
here too (Odifreddi tells   ``{\sl
[Gerald] Edelman  quotes Turing as a precursor of his work 
on morphogenesis\/}").

What is going on here is fascinating. What 
is taking shape seems to be the 1936 paradigm shift renewed. In the 
1920s and 1930s we saw a fracturing of the comfortable 
picture of how science could bring predictability to 
a complex universe, unaccompanied by any overall concept 
of the underlying mathematical structures. 
 We now see a  
coming together (at times faltering, at times confused) of science
and mathematics to replace the Laplacian  model of science with
one whose complexities  match those of the real world. At the
root of this  is both the technical legacy of Turing, and the 
kind of unified approach to scientific problems that 
is so characteristic of Turing's own thinking. 

The role of the classical computability theorist here 
is, as yet, peripheral. There is little awareness and
understanding  of much recent research directly
descended  from Turing's work on incomputability and his notion of 
relative computability, beyond the tiny community of specialists.  
(Turing's computability-theoretic heirs are significantly absent 
from Teuscher's recent book \cite{teuscher} ---  otherwise 
such a useful window on this Turing rennaissance.) 
Much of the work being done has a 
mathematically naive aspect on which sceptics such as Martin Davis
have  been quick to pounce. 
But the separation since Turing --- during   
  the recursion theoretic era  --- between computability theory 
and its real-world counterpart, is being slowly repaired, 
as is a lingering distrust of current outcomes, on both sides of the divide. 

\section{Turing's final thoughts on the mind as machine?} 

So Alan Turing turns out to be a pivotal figure in 
 a scientific revolution which has 
not yet run its course. And right to the end he 
was grappling with the complexities of 
physical computability and its mathematical models. 
For Turing the scientist, the human mind was his laboratory. 
And, as I have suggested, his 1939 paper provided the theoretical framework 
for all his subsequent thinking on the mind as 
computer. But a measure of the extent to which he was 
prepared to engage with everyday problems is the 
way he was able to communicate his science 
to everyday people. Here is the final paragraph of 
Alan Turing's article \cite{turing54} on  
 {\em Solvable and Unsolvable Problems}, in the popular 
  {\em Penguin Science 
News\/} back in his final year of 1954, (number 31, p.23):

{\rightskip=14pt \leftskip=14pt \hskip 0mm {\sl 
``The results which have been described in this 
article are mainly of  a negative character, setting 
certain bounds to what we can hope to achieve purely 
by reasoning. These, and some other results of mathematical 
logic may be regarded as going some way towards 
a demonstration, within mathematics itself, 
of the inadequacy of `reason' \underbar{unsu}pp\underbar{orted 
b}y\underbar{ common 
sense}." \hskip 3mm{\rm (my underlining)}}\vskip 1mm} 

After my talk on which this article is based, 
Jack Copeland asked me whether I thought that Turing meant 
 ``intuition"  when he referred here to 
 ``common sense" --- ``if so, I agree". Yes, despite Turing's 
characteristically cautious phraseology, the above quotation 
clearly comes from the inventor and 
interpreter of ordinal logics.  Certainly not someone 
who was content with any existing algorithmic model 
of the human mind. 

And in his final months he was renewing his interest in 
the parallel mysteries of quantum theory, as indicated 
by his cryptic \textit{Messages from the Unseen World\/} postcards 
 to Robin Gandy (see
\cite{hodges}, p.512). Alan was  still the writer to Mrs Morcom of those youthful
speculations  on the links between incomputability in the mind and at 
the quantum level. Here is Robin Gandy in a letter   
 to Max Newman in June 1954:
 
{\rightskip=14pt \leftskip=14pt \hskip 0mm {\sl 
``During this spring [Turing] spent some time inventing 
a new quantum mechanics \dots he produced a slogan `Description 
must be non-linear, prediction must be linear'."}
\end{document}